\documentclass[12pt]{amsart}
\usepackage{amsfonts,latexsym,amsthm,amssymb,graphicx}
\usepackage[all]{xy}

\usepackage[english,francais]{babel}
\newtheorem{theorem}{Theorem}[section]

\newtheorem{lemma}[theorem]{Lemme}
\newtheorem{thrm}[theorem]{Th\'{e}or\`{e}me}

{\theoremstyle{definition} \newtheorem{defin}[theorem]{D\'efinition}}
{\theoremstyle{remark} 
\newtheorem{example}[theorem]{Exemple}}

\newenvironment{demo}{\noindent \textsl{Preuve.}}{\qed}

\newcommand{\Cbb}{{\mathbb{C}}}

\newcommand{\one}{1\hskip-3.5pt1}
\newcommand{\cst}{c_*}
\newcommand{\ctil}{{\tilde c}}

\title[Co\"incidence des classes de Schwartz et de MacPherson]{Une nouvelle 
preuve de la co\"incidence des classes d\'efinies par M.-H.~Schwartz et par R.~MacPherson}
\author{Paolo Aluffi et Jean-Paul Brasselet}
\address{
Mathematics Department, 
Florida State University,
Tallahassee FL 32306, U.S.A.
}
\email{aluffi@math.fsu.edu}
\address{
IML -- CNRS, 
Case 907, Luminy, 
13288 Marseille Cedex 9, France
}
\email{jpb@iml.univ-mrs.fr}

\begin{document}

\begin{abstract}
Nous donnons une courte d\'emonstration de ce que les classes des vari\'et\'es 
singuli\`eres d\'efinies par Marie-H\'el\`ene Schwartz au moyen des ``champs radiaux" 
co\"incident avec la notion fonctorielle d\'efinie par Robert MacPherson.
\end{abstract}

\maketitle

\begin{flushright}
{\sl A Marie-H\'el\`ene Schwartz}
\end{flushright}


\newcommand{\oS}{{\overline S}}
\newcommand{\oT}{{\overline T}}
\newcommand{\cE}{{\mathcal E}}

\section{Introduction}\label{intro}
Marie-H\'el\`ene Schwartz (\cite{MR35:3707}, \cite{MR32:1727}) 
a d\'efini, au milieu des ann\'ees 60,  une notion de {\em  classes de Chern 
pour les vari\'et\'es singuli\`eres\/} en 
cohomologie relative, ou, via dualit\'e d'Alexander,  en homologie \`a coefficients
entiers. La m\'ethode utilis\'ee est  la th\'eorie d'obstruction appliqu\'ee \`a des 
champs de rep\`eres  
dont le comportement est control\'e le long de la partie singuli\`ere.

Par la suite, le travail d'Alexander Grothendieck en vue de montrer un
`th\'eor\`eme de Riemann-Roch discret' a \'et\'e \`a l'origine de
la conjecture de l'existence (en caract\'eristique z\'ero) d'une
th\'eorie fonctorielle de classes de Chern,  comme
transformation naturelle du foncteur des fonctions constructibles
dans une th\'eorie d'homologie convenable. Cette conjecture est
connue sous le nom de conjecture de Deligne-Grothendieck.
Celle-ci a \'et\'e r\'esolue, dans les ann\'ees 70,
par Robert MacPherson \cite{MR50:13587}, lequel a ainsi donn\'e
une autre construction de classes de Chern pour les
vari\'et\'es singuli\`eres.

Le travail de R. MacPherson  est ind\'ependant de celui de M.H. Schwartz, 
bien que les deux notions aient des points communs : d'une part, les deux notions de classes 
sp\'ecialisent aux classes de Chern usuelles pour les vari\'et\'es lisses ; d'autre part, le fait que 
le r\'esultat de M.H. Schwartz \'etende le Th\'eor\`eme de Poincar\'e-Hopf aux 
vari\'et\'es singuli\`eres peut \^etre consid\'er\'e comme un aspect de la propri\'et\'e de 
fonctorialit\'e satisfaite par les classes de MacPherson.

Il \'etait donc naturel de conjecturer que ces deux notions co\"incident, ce qui a \'et\'e 
d\'emontr\'e en 1979 par M.H. Schwartz et le second auteur du pr\'esent article :

\begin{thrm}[\cite{MR83h:32011}]\label{main}
Les classes de Chern homologiques de Schwartz et de  MacPherson sont \'egales.
\end{thrm}

La m\'ethode de 
d\'emonstration consiste \`a relier les indices des champs radiaux et les classes 
de Schwartz aux ingr\'edients principaux de la construction des classes de MacPherson, 
\`a savoir l'{\em  obstruction d'Euler locale\/} et {\em la classe de Chern-Mather.\/}

Dans cet article, nous montrons plus directement l'\'egalit\'e des classes de 
Schwartz et de MacPherson, sans nous servir des autres invariants des
singularit\'es. Ce nouveau point de vue est inspir\'e par  une nouvelle 
expression de la notion de fonctorialit\'e des classes de Chern en termes de classes 
d\'efinies pour les vari\'et\'es lisses (mais non n\'ecessairement 
{\em compl\`etes\/}), obtenue par le premier auteur  \cite{math.AG/0507029}. 
En fait, des d\'eveloppements ult\'erieurs ont rendu l'argument ind\'ependant 
de cette r\'ef\'erence, ce qui permet de pr\'esenter la pr\'esente version 
de fa\c con autonome.

Dans la section \ref{crite}, on montre qu'une classe 
d\'efinie pour les vari\'et\'es $X$ (\'eventuellement singuli\`eres) co\"incide 
n\'ecessairement avec la notion fonctorielle si
\begin{itemize}
\item elle peut \^etre \'ecrite comme somme des contributions des strates  
d'une stratification de~$X$;
\item la contribution d'une strate non singuli\`ere $S$ est pr\'eserv\'ee par 
morphismes respectant les stratifications ; et 
\item la classe co\"\i ncide avec la classe de Chern du fibr\'e tangent d'une 
vari\'et\'e compl\`ete non singuli\`ere.

\end{itemize}

Dans la section \ref{schwartz} on remarque que les classes d\'efinies par Marie-H\'el\`ene 
Schwartz satisfont ces conditions, et le fait qu'elles co\"\i ncident avec les classes 
de MacPherson en r\'esulte imm\'ediatement.


\section{Une caract\'erisation des classes de Chern fonctorielles}
\label{crite}

Dans cette section, nous travaillons sur un corps arbitraire alg\'ebriquement clos 
de caract\'eristique z\'ero, et dans le groupe de Chow $A_*$.  Les r\'esultats sont valables 
{\em a fortiori\/} pour les vari\'et\'es complexes, en homologie.

\subsection{}
Soit $\ctil(X)\in A_*X$ une classe d\'efinie pour toutes les vari\'et\'es (\'eventuellement 
singuli\`eres) $X$, comme une somme de contributions d'une d\'ecomposition de $X$ 
en union finie disjointe de vari\'et\'es (\'eventuellement incompl\`etes) 
{\em non singuli\`eres\/} $S_i$ :
$$X=\amalg_i S_i \leadsto \ctil(X)=\sum_i \ctil(S_i,X)\quad;$$
on dit que de telles d\'ecompositions sont {\em admissibles\/} (pour $\ctil$).
Comme nous le verrons, dans certains cas {\em toute\/} d\'ecomposition de $X$ 
comme union finie disjointe de sous-vari\'et\'es non singuli\`eres peut \^etre 
admissible. Dans d'autres situations, des restrictions 
peuvent \^etre impos\'ees sur les d\'ecompositions pour \^etre admissibles : 
par exemple, on peut demander aux vari\'et\'es $S_i$ d'\^etre \'el\'ements d'une  
stratification de Whitney de $X$.

Nous supposons que les strates d'un diviseur \`a croisements normaux forment une 
d\'ecomposition admissible. Plus pr\'ecis\'ement : si $D$ est un diviseur \`a croisements 
normaux simples, de composantes non singuli\`eres $D_j$, $j\in J$, d'une vari\'et\'e 
non singuli\`ere $Y$, on suppose que 
$$\amalg_{I\subset J} D_I^\circ$$
est une d\'ecomposition admissible de $Y$, o\`u $D_I^\circ$ d\'esigne 
$$(\cap_{j\in I} D_j)\smallsetminus (\cup_{j\not\in I} D_j)\quad.$$
(Par exemple, $D_\emptyset^\circ$ est le compl\'ementaire de $D$ dans $Y$).

Cette propri\'et\'e sera imm\'ediatement satisfaite pour les d\'ecompositions que nous consid\'ererons.

\begin{defin}
On dit que la classe $\ctil(X)$ est {\em localement d\'etermin\'ee\/} si 
la condition suivante est r\'ealis\'ee :
\begin{itemize}
\item Si $f: Y \to X$ est un morphisme propre, $S$ et $T:=f^{-1}(S)$ 
resp.~sont \'el\'ements d'une d\'ecomposition admissible de $X$, resp. $Y$, 
et $f$ se restreint \`a un isomorphisme $T \to S$, alors 
$$f_* \ctil(T,Y)=\ctil(S,X)\quad.$$
\end{itemize}
\end{defin}

\begin{example}[La  classe fonctorielle]
On d\'esigne par $\cst(X)\in A_*(X)$ la classe d\'efinie par MacPherson dans
\cite{MR50:13587} (voir \cite{MR85k:14004}, \S19.1.7, 
pour l'adaptation de la d\'efinition au groupe de Chow $A_*(X)$, et 
\cite{MR91h:14010} pour l'extension aux corps arbitraires alg\'e\-bri\-que\-ment clos
de caract\'eristique z\'ero). Rappelons que cette classe est la valeur 
$$\cst(X):=c_*(\one_X)\in A_*X$$
prise, sur la fonction caract\'eristique constante $\one_X$, par une transformation naturelle 
$$c_*: F \leadsto A_*$$
du foncteur des fonctions constructibles dans le foncteur groupe de  Chow. Ici 
$F(X)$ est le groupe des fonctions constructibles \`a valeurs enti\`eres sur $X$ ;
si $g:Y \to X$ est une application propre, 
l'image directe 
$g_*(\varphi)$ de la 
fonction constructible $\varphi=\sum_Z m_Z \one_Z \in F(Y)$ est d\'efinie comme  
la fonction sur $X$ dont la valeur en $p\in X$ est 
$$g_*(\varphi)(p) := \sum m_Z \, \chi(g^{-1}(p)\cap Z)\quad.$$
Ici $\chi$ d\'esigne la caract\'eristique d'Euler-Poincar\'e  topologique, sur $\Cbb$ ; pour 
l'extension aux autres corps de caract\'eristique z\'ero voir \cite{MR91h:14010}
ou \cite{math.AG/0507029}.

Si $S$ d\'esigne une 
sous-vari\'et\'e 
quelconque (en particulier non singuli\`ere) de $X$, on pose 
$$\cst(S,X):=c_*(\one_S)\in A_*X\quad.$$

Alors {\em toute\/} d\'ecomposition de $X$ comme union finie disjointe de 
sous-vari\'et\'es  
non singuli\`eres est admissible pour $\cst$.
En effet, si $X=\amalg_i S_i$ alors $\one_X=\sum_i \one_{S_i}$ et donc
$$\cst(X)=c_*(\one_X)=\sum_i c_*(\one_{S_i})=\sum_i \cst(S_i,X)\quad.$$

De plus, $\cst$ est localement d\'etermin\'ee. En effet,  soit $f:Y \to X$ une application
propre qui se restreint \`a un isomorphisme $T \to S$. Alors
$$f_*\cst(T,Y)=f_*c_*(\one_T)=c_* f_*(\one_T)=c_*(\one_S)=\cst(S,X)$$
puisque $c_*$ est une transformation naturelle, et par d\'efinition de l'image 
directe des
fonctions constructibles.
\end{example}

\begin{thrm}\label{crit}
Supposons que $\ctil$ soit localement d\'etermin\'ee, et que 
$$\ctil(V)=c(TV)\cap [V]$$
pour toute vari\'et\'e non-singuli\`ere {\em compl\`ete\/} $V$. Alors $\ctil$
co\"\i ncide avec la classe fonctorielle $\cst$.
\end{thrm}

\begin{demo}
Soit $X=\amalg_i S_i$ une d\'ecomposition admissible de $X$ pour $\ctil$;
il suffit de montrer que $\ctil(S_i,X)=\cst(S_i,X)$. 

Si $S$ est un \'el\'ement de la d\'ecomposition, notons $\oS$ son adh\'erence dans $X$, et consid\'erons 
$f:Y \to \oS$ une r\'esolution plong\'ee de $\oS$, telle que le compl\'ement 
de $T:=f^{-1}(S)$ dans $Y=\oT$ soit un diviseur $D$ \`a croisements normaux, de  composantes 
non singuli\`eres  $D_j$, $j\in J$.
Comme $\ctil$ et $\cst$ sont toutes deux localement d\'etermin\'ees, il suffit de montrer que 
$$\ctil(T,Y)=\cst(T,Y)\in A_*Y\quad.$$
On a
$$Y=\amalg_{I\subset J} D_I^\circ\quad,$$
d'o\`u  
$$\ctil(Y)=\sum_{I\subset J} \ctil(D_I^\circ, Y)\quad,$$
ce qui, par le principe d'``inclusion-exclusion",
donne 
$$\ctil(T,Y)=\ctil(D_\emptyset^\circ,Y)=\sum_{I\subset J}(-1)^{|I|}\ctil(D_I,Y)
\quad,$$
o\`u $D_I=\cap_{j\in I} D_j$, et nous notons par $\ctil(D_I,Y)$ la somme 
$\sum_{K\supset I}\ctil(D_K^\circ,Y)$. 
Chaque $D_I$ est compl\`ete et non singuli\`ere 
(puisque $D$ est un diviseur \`a croisements normaux), donc l'inclusion
$\iota:D_I \to Y$ est propre. Pour $K\supset I$ on a
$$\ctil(D_K^\circ,Y)=\iota_* \ctil(D_K^\circ,D_I)\quad,$$
puisque $\ctil$ est localement d\'etermin\'ee, et donc
$$\ctil(D_I,Y)=\sum_{K\supset I}\iota_* \ctil(D_K^\circ, D_I)=\iota_*\ctil(D_I)
=\iota_* (c(T D_I)\cap [D_I])\quad.$$
En utilisant les m\^emes arguments, on a la m\^eme expression pour 
$\cst(T,Y)$, d'o\`u le r\'esultat.
\end{demo}

 
\section{Les classes de Schwartz}\label{schwartz}
 
Nous supposons maintenant que le corps de base est $\Cbb$, 
et nous travaillons en homologie \`a coefficients entiers ; les consid\'erations de la section 
pr\'ec\'edante s'appliquent dans ce contexte, via l'application canonique du groupe
de Chow dans l'homologie.
 
\subsection{}
Notons $\ctil$ la classe d\'efinie par Marie-H\'el\`ene Schwartz
dans \cite{MR35:3707}, \cite{MR32:1727}. Cette classe co\"\i ncide avec la classe 
de Chern (homologique) du fibr\'e tangent pour les vari\'et\'es non singuli\`eres compl\`etes,
et se calcule, dans le cas g\'en\'eral, comme somme des contributions des strates 
$S_i$ d'une stratification de Whitney de $X$:
$$\ctil(X):=\sum_i \ctil(S_i,X)\quad.$$
Autrement dit, les stratifications de Whitney sont des d\'ecompositions admissibles pour $\ctil$.

Rappelons bri\`evement la d\'efinition des classes de Schwartz. La vari\'et\'e $X$
est plong\'ee dans une vari\'et\'e non singuli\`ere $M$, stratifi\'ee par $M\setminus X$
et par les strates $S_i$ de $X$. Pour $x\in M$, on note $E(x)$ le sous-espace
de l'espace tangent $T_xM$ form\'e des vecteurs tangents
\`a la strate de $M$ contenant $x$. La collection des sous-espaces 
$E(x)$ d\'etermine un sous-espace $E$ du fibr\'e tangent 
$TM$. Ce n'est plus un fibr\'e mais la  notion d'une section de $E$ est bien d\'efinie,
comme section de $TM$ \`a valeurs dans $E$.
On consid\`ere alors l'espace $E_r$ des $r$-rep\`eres ordonn\'es de 
vecteurs tangents aux strates de $M$. C'est un sous-espace du fibr\'e 
$T_rM$ des $r$-rep\`eres ordonn\'es de $TM$.

Consid\`erons aussi une triangulation $K$ de $M$ compatible avec la 
stratification, et une d\'ecomposition cellulaire  $D$ de $M$, duale de $K$. 
Les cellules de $D$ sont transverses aux strates.

En dimension $r-1$, la classe de  Schwartz est d\'efinie comme l'obstruction 
\`a la construction d'une section particuli\`ere $Z_r$ de $E_r$ sur $X$, appel\'ee {\em champ  
de $r$-rep\`eres radial}. C'est en particulier une section de 
$T_rM$ sans singularit\'e sur le squelette de dimension (r\'eelle) 
$2m-2r+1$ de la d\'ecomposition cellulaire duale 
$D$ et avec singularit\'es isol\'ees sur le $D$--squelette de dimension $2(m-r+1)$.
La construction est faite par r\'ecurrence sur la dimension des strates : 

Les strates de dimension complexe (strictement) inf\'erieure \`a $r-1$ n'apparaissent pas. 
Consid\'erons une strate de dimension $(r-1)$. Les cellules de dimension (r\'eelle) 
$2(m-r+1)$ intersectent une telle strate 
en un point (lorsque l'intersection est non  vide). Sur une telle cellule, le $r$-rep\`ere 
est construit comme champ radial, et donc d'indice $+1$. 
Si $\dim_{\Cbb} S = r-1$, on a alors : 
$$\ctil(S,X) = \sum_{K_\alpha\subset S} K_\alpha^{r-1}$$

Supposons que la construction soit d\'ej\`a effectu\'ee sur les strates de 
dimension 
inf\'erieure
\`a celle de $S$. Le $r$-rep\`ere 
$Z_r$ est d\'efini sur les strates du bord de $S$, avec des singularit\'es 
$a_j$. Notons  $U$ le voisinage ``cellulaire" de $\partial S$, 
compos\'e des cellules de $D$ duales des simplexes de $\partial S$.
La construction s'effectue de fa\c con \`a ce que, dans $U$, l'extension de 
$Z_r$, encore not\'ee $Z_r$ n'a pas d'autre singularit\'e que les points $a_j$. 
D'apr\`es la construction des champs radiaux et comme prouv\'e dans 
 Schwartz \cite{MHS}, le $r$-rep\`ere  $Z_r$ satisfait aux propri\'et\'es suivantes :

1) L'indice $I(Z_r,a_j)$ est le m\^eme, calcul\'e comme section du fibr\'e tangent 
\`a la strate de $a_j$ ou comme section de $TM$, 

2) Le rep\`ere $Z_r$ est sortant  du voisinage cellulaire $U$. 
Cela signifie que sur $\partial U \cap S$, 
le $r$-rep\`ere $Z_r$ est tangent \`a $S$ et entrant dans 
$S\setminus U$, 

3) Deux champs radiaux sont homotopes sur $\partial U \cap S$, 
comme sections de $T_rS$. 

On \'etend alors $Z_r$,  d\'efini dans $S\cap U$, dans l'int\'erieur de $S$,   
avec singularit\'es isol\'ees $a$ situ\'ees dans l'intersection de $S$ 
avec les $D$-cellules 
de dimension $2(m-r+1)$.
Autrement dit, les points $a$ sont  situ\'es 
dans l'intersection de $S$ avec les $D$-cellules 
duales des $K$-simplexes de dimension $2(r-1)$ contenus dans $S\setminus U$. 

La contribution $\ctil(S,X)$ de la strate $S$ se calcule en termes des indices 
$I(Z_r,a)$ en ces points singuliers $a$ de $Z_r$, situ\'es dans $S\setminus \partial S$. 
$$\ctil(S,X):=\sum_{K_\alpha\subset S\setminus \partial S} \mu_\alpha
K_\alpha^{r-1}$$
o\`u les simplexes $K_\alpha$ sont de dimension $2(r-1)$ et
$$\mu_\alpha=\sum_{a\in D_\alpha\cap S} I(Z_r,a)\quad,$$
o\`u la cellule $D_\alpha$ est duale de $K_\alpha$.

Une cons\'equence des Propri\'et\'es 2) et 3) est le Lemme suivant :

\begin{lemma}\label{indepdt}
La contribution 
$\ctil(S,X)$ ne d\'epend pas de la  construction  du rep\`ere  radial
sur les strates du bord $\partial S$. Autrement dit,  
deux champs radiaux 
$Z_r$ et $Z'_r$ d\'efinis sur  $S\cap \partial U$ donnent lieu \`a des cycles  \'equivalents
$$\sum_{K_\alpha\subset S} \mu_\alpha K_\alpha^{r-1}
\sim \sum_{K_\alpha\subset S} \mu'_\alpha K_\alpha^{r-1}$$
\end{lemma}

Remarquons que, via l'isomorphisme d'Alexander $H_{2(r-1)}(S) \cong H^{2p}(M, M\setminus S)$, 
avec $p=m-r+1$, la classe $\ctil(S,X)$ correspond \`a la classe not\'ee $\hat c^{2p} (S)$ 
dans Schwartz. 

D'apr\`es le Th\'eor\`eme~\ref{crit}, pour montrer le Th\'eor\`eme~\ref{main} il 
suffit de montrer que les classes de Chern de Schwartz sont localement d\'etermin\'ees,
ce qui est l'objet du Lemme suivant : 

\begin{lemma}
Si $f: Y \to X$ est un morphisme propre, $S$ et $T:=f^{-1}(S)$ 
resp.~sont \'el\'ements de d\'ecompositions admissibles de $X$, $Y$ resp., 
et $f$ se restreint \`a un isomorphisme $T \to S$,  alors
$$f_* \ctil(T,Y)=\ctil(S,X)\quad.$$
\end{lemma}

\begin{demo}
Notons $\partial S$ et $\partial T$ les bords de 
$S$ et $T$ respectivement, unions de strates de $X$ et $Y$
et notons $k$ la dimension (complexe) commune de $S$ et $T$.
Supposons effectu\'ee la construction des classes de Schwartz de $X$ et $Y$ 
par r\'ecurrence sur la dimension des strates jusqu'en dimension $k-1$ incluse. 
Cette construction donne lieu \`a 
deux champs de $r$-rep\`eres :  l'un $Z_r$, d\'efini sur un voisinage cellulaire $U$ de 
$\partial S$, sortant de $U$, et l'autre $W_r$, d\'efini  sur un voisinage cellulaire
$V$ de $\partial T$, sortant de $V$. Nous avons donc la situation suivante :

Sur $S\cap \partial U$ le $r$-rep\`ere $Z_r$ est tangent \`a $S$ et entrant dans $S\setminus U$.  
Sur $T\cap \partial V$ le $r$-rep\`ere  $W_r$ est tangent \`a $T$ et entrant dans $T\setminus V$.  
D'apr\`es l'hypoth\`ese du Th\'eor\`eme, on a un isomorphisme des paires 
$(S\setminus U, S\cap \partial U) 
\cong (T\setminus V, T\cap \partial V)$. On conclut par le 
Lemme~\ref{indepdt}.
\end{demo}
 


\end{document}